\documentclass[12pt]{amsart}
\font\emailfont=cmtt10

\usepackage{amsmath,amsthm,amsfonts,amscd,flafter,epsf}

\newcommand{\rk}{\mathrm{rk}}
\newcommand{\HF}{HF}

\newtheorem{theorem}{Theorem}[section]
\newtheorem{prop}[theorem]{Proposition}
\newtheorem{cor}[theorem]{Corollary}

\newtheorem{lemma}[theorem]{Lemma}

\def\endproof{\relax\ifmmode\expandafter\endproofmath\else
  \unskip\nobreak\hfil\penalty50\hskip.75em\hbox{}\nobreak\hfil\bull
  {\parfillskip=0pt \finalhyphendemerits=0 \bigbreak}\fi}
\def\endproofmath$${\eqno\bull$$\bigbreak}
\def\bull{\vbox{\hrule\hbox{\vrule\kern3pt\vbox{\kern6pt}\kern3pt\vrule}\hrule}}

\newcommand{\Q}{\mathbb{Q}}

\newcommand{\Z}{\mathbb{Z}}

\newcommand{\Zmod}[1]{\Z/{#1}\Z}

\newcommand{\Ker}{\mathrm{Ker}}

\newcommand{\cm}{\cdot}

\newcommand{\ModSWfour}{\mathcal{M}}
\newcommand{\ModFlow}{\ModSWfour}

\newcommand{\SpinC}{{\mathrm{Spin}}^c}

\newcommand\abuts\Rightarrow
\newcommand\Sym{\mathrm{Sym}}

\newcommand\HFpRed{\HFp_{\red}}

\newcommand\relspinc{\underline{\spinc}}

\newcommand\ModSphere{\ModFlow\left({\mathbb S}\longrightarrow 
\Sym^{g-1}(\Sigma_{1})\times \Sym^2(\Sigma_{2})\right)}
\newcommand\ModSpheres\ModSphere

\newcommand\CFp{\CFb}

\newcommand\HFred{\HF_{\rm red}}

\newcommand{\red}{\mathrm{red}}

\newcommand\HFp{\HFb}

\newcommand\CFb{CF^+}

\newcommand\HFb{HF^+}

\newcommand\UnparModSp{\widehat \ModSp}
\newcommand\UnparModFlow\UnparModSp
\newcommand\Mod\ModSp

\newcommand{\spinc}{\mathfrak s}

\newcommand{\spinct}{\mathfrak t}

\newcommand\ModMaps{\mathcal M}
\newcommand\ModSp\ModMaps

\newcommand\Dual{\mathcal D}
\newcommand\Duality\Dual

\newcommand\Hred{H^{\red}}
\newcommand\InjMod{{\mathcal T}^+}

\newcommand\Dp{D^+}
\newcommand\Xp{{\mathbb X}^+}
\newcommand\BigAp{{\mathbb A}^+}
\newcommand\BigBp{{\mathbb B}^+}
\newcommand\Ap{A^+}
\newcommand\Bp{B^+}
\newcommand\vertp{v^+}
\newcommand\horp{h^+}

\newcommand\spincrel\relspinc

\newcommand\CFK{\mathrm{CFK}}
\newcommand\HFK{\mathrm{HFK}}

\newcommand\CFKinf{\CFK^{\infty}}

\newcommand\HFKa{\widehat\HFK}

\keywords{Floer homology, trefoil knot, figure eight knot, Dehn surgery}

\title[{Trefoil and figure eight knot}] 
{The Dehn surgery characterization of the trefoil
and the figure eight knot}

\author[Peter Ozsv{\'a}th]{Peter Ozsv\'ath}
\address{Department of
Mathematics, Columbia University, 
New York, NY 10027 \newline
\indent{\emailfont{petero@math.columbia.edu}}}
\thanks{PSO was supported by NSF grant number DMS-050581}

\author[Zolt{\'a}n Szab{\'o}]{Zolt{\'a}n Szab{\'o}} 
\address{Department of
Mathematics, Princeton University, New Jersey 08544 \newline
\indent{\emailfont{szabo@math.princeton.edu}}}
\thanks{ZSz was supported by NSF grant number DMS-0406155}

\begin{document}

\begin{abstract}  
  We give a Dehn surgery characterization of the trefoil and the
  figure eight knots. These results are gotten by combining surgery
  formulas in Heegaard Floer homology from an earlier paper with the
  characterization of these knots in terms of their knot Floer
  homology given in a recent paper of Ghiggini.
\end{abstract}

\maketitle
\section{Introduction}

In~\cite{KMOSz}, it was shown that the unknot is characterized
by its Dehn surgeries, in the sense that if Dehn surgery with some slope
on a knot $K$ in $S^3$ is orientation-preserving homeomorphic to Dehn surgery 
with the same slope on the unknot, then $K$ is in fact unknotted.

A proof of this fact can be given using Heegaard Floer
homology~\cite{HolDisk}.  Specifically, there is a Heegaard Floer
homology theory for knots introduced in~\cite{Knots}
and~\cite{RasmussenThesis}. Surgery formulas for this
invariant~\cite{RatSurg} allow one to express the Heegaard Floer
homology for $p/q$-Dehn surgery of $K$ in terms of this knot Floer
homology of $K$. The hypothesis that $p/q$ surgery on $K$ agrees with
that of the unknot forces the knot Floer homology of $K$ to agree with
that of the unknot. Combining this with the fact that knot Floer
homology detects the unknot~\cite{GenusBounds}, the Dehn surgery
characterization of the unknot follows.

In a beautiful recent paper, Ghiggini~\cite{Ghiggini} shows that
Heegaard Floer homology also detects the trefoil and the figure eight
knot. Appealing to the same strategy outlined above, in the form of
the surgery formulas for knot Floer homology, we obtain here a similar
Dehn surgery characterization of both of these knots. Specifically, we
have the following:

\begin{theorem}
  \label{thm:Trefoil}
  Let $T$ be a trefoil knot.  If $K$ is a knot with the property that
  there is a rational number $r$ and an orientation-preserving
  diffeomorphism $S^3_r(K)\cong S^3_r(T)$, then $K$ is in fact the
  trefoil $T$.
\end{theorem}

In a similar vein, we have the following 

\begin{theorem}
  \label{thm:FigEight}
  Let $S$ be the figure eight knot.  If $K$ is a knot with the
  property that there is a rational number $r$ and an
  orientation-preserving diffeomorphism $S^3_r(K)\cong S^3_r(S)$, then
  $K$ is in fact the figure eight knot $S$.
\end{theorem}

The condition that the diffeomorphism preserves orientations is
important, here. For example, there are identifications
$S^3_{+1}(T_\ell)\cong -S^3_{+1}(S)$, and also $S^3_{+5}(T_r)\cong
-S^3_{+5}(O)$, where here $O$ is the unknot.

Similarly, the condition that the surgery coefficient is fixed is also
crucial; $S^3_{1/n}(K_0)$ where $K_0$ is a trefoil or the 
figure eight knot can be realized alternatively as $+1$ surgery on 
a suitable twist knot.

In his paper, Ghiggini proves that the trefoil is the only knot in
$S^3$ which admits a surgery giving the Poincar\'e homology sphere.  A
consequence of Theorems~\ref{thm:Trefoil} and~\ref{thm:FigEight}, we
obtain a similar result for the Brieskorn sphere $\Sigma(2,3,7)$.

\begin{cor}
\label{cor:S237}
The only surgeries on knots in $S^3$ which realize the Brieskorn
sphere  $\Sigma(2,3,7)$ (with either orientation) are $S^3_{-1}(T_r)\cong S^3_{+1}(S)\cong \Sigma(2,3,7)$
and 
$S^3_{+1}(T_\ell)\cong S^3_{-1}(S)\cong -\Sigma(2,3,7)$. 
\end{cor}

In Section~\ref{sec:Background}, we review the relevant aspects of
Heegaard Floer homology which are used in the proofs of the above
results. In Section~\ref{sec:Proof} we give the proofs of the above two 
theorems and the corollary.

\section{Background}
\label{sec:Background}

\subsection{Heegaard Floer homology}

In its most basic form, Heegaard Floer homology is a $\Zmod{2}$-graded
Abelian group associated to a three-manifold, but it comes in several variants
and can be endowed with additional structure~\cite{HolDisk}.

In this paper, we will consider primarily the version
$\HFp(Y)$ for rational homology three-spheres $Y$. This group admits
a splitting according to $\SpinC$ structures over $Y$
$$\HFp(Y)\cong \bigoplus_{\spinct\in\SpinC(Y)} \HFp(Y,\spinct).$$
Moreover, $\HFp(Y,\spinct)$ is equipped with an absolute grading
(defined in~\cite[Section~\ref{HolDiskFour:sec:AbsGrade}]{HolDiskFour}
and studied extensively in~\cite{AbsGraded}). Recall that there is a natural
involution on the space of $\SpinC$ structures over $Y$, denoted
$\spinct\mapsto{\overline \spinct}$. There is a corresponding isomorphism
\begin{equation}
\label{eq:Conjugation}
\HFp(Y,\spinct)\cong \HFp(Y,{\overline \spinct}),
\end{equation}
which, in the case of rational homology spheres $Y$, is an isomorphism
of $\Q$-graded Abelian groups.

The group $\HFp(Y,\spinct)$ has the following algebraic structure.
$$\HFp(Y,\spinct)=\bigoplus_{d\in\Q} \HFp(Y,\spinct).$$ In fact,
$\HFp(Y,\spinct)$ is supported only in rational degrees $d$ within
some fixed equivalence class (depending on $Y$ and $\spinct$) modulo
the integers.  For any degree $d\in\Q$, $\HFp_d(Y,\spinct)$ is a
finitely generated $\Z$-module. Moreover, $\HFp(Y,\spinct)$ is
endowed with an endomorphism $U$ which lowers degree by
$2$, i.e.
$$U\colon \HFp_{d}(Y,\spinct)\longrightarrow \HFp_{d-2}(Y,\spinct).$$ 
Moreover, $\HFp_{d}(Y,\spinct)=0$ for all sufficiently small $d$. 
Finally, for any sufficiently large rational number $d_0$,
if we consider the quotient module $\HFp_{\geq d_0}(Y,\spinct)$
of $\HFp(Y,\spinct)$
generated by all elements with degree greater than $d_0$, then 
that module is isomorphic to the $\Z[U]$-module
$$\InjMod=\frac{\Z[U,U^{-1}]}{U \cm \Z[U]}.$$

From the above properties, it is clear that there is a canonical short exact sequence
$$\begin{CD}
0@>>>\InjMod@>>>\HFp(Y,\spinct)@>>>\HFpRed(Y,\spinct)@>>> 0,
\end{CD}$$
where here $\HFpRed(Y,\spinct)$ is $\Z[U]$ module which is a finitely
generated $\Z$-module. Moreover, we obtain a three-manifold invariant,
$$d\colon \SpinC(Y)\longrightarrow \Q,$$ the {\em correction terms} of
$Y$, where $d(Y,\spinct)$ is the minimal $\Q$-grading of any
homogeneous element of $\HFp(Y,\spinct)$ in the image of $\InjMod$.

The Floer homology group $\HFp(Y,\spinct)$ also inherits a $\Zmod{2}$-grading;
a non-zero element in  $\HFp_d(Y,\spinct)$ has {\em even parity} if $d\equiv d(Y,\spinct)\pmod{2\Z}$,
and it has {\em odd parity} if $d\equiv d(Y,\spinct)+1\pmod{2\Z}$.
The correction terms are analogous to a gauge-theoretic invariant
introduced by Fr{\o}yshov~\cite{Froyshov}; for more information on the
correction terms, see~\cite{AbsGraded}.

If $C$ is a chain complex of $\Z[U]$-modules, we can 
form 
$$\Hred(C) = \lim_{\stackrel{\leftarrow}{d}} \frac{H_*(C)}{U^d H_*(C)}.$$
Then, if $\CFp(Y,\spinct)$ is the chain complex calculating $\HFp(Y,\spinct)$, then
$\Hred(\CFp(Y,\spinct))=\HFpRed(Y,\spinct)$.

\subsection{Knot Floer homology and the surgery formula}
Heegaard Floer homology can be extended as in ~\cite{Knots}
and~\cite{RasmussenThesis} to invariants for null-homologous knots $K$
in closed three-manifolds.  We restrict attention to the case where
the ambient three-manifold is the three-sphere $S^3$. We recall now
the notation of knot Floer homology, following~\cite{Knots}. This data
can be used to calculate Heegaard Floer homology groups of Dehn
fillings of $S^3$ along $K$, cf.~\cite{RatSurg}. After setting up
notation for knot Floer homology, we state this surgery formula.

A $\Z\oplus\Z$-filtered chain complex is a free Abelian group which
splits as a direct sum $C=\bigoplus_{(i,j)\in \Z\oplus Z} C\{(i,j)\}$
and which is endowed with a boundary operator which carries elements
in $C\{(i,j)\}$ to elements in
$$\bigoplus_{(i',j')\leq (i,j)}C\{(i',j')\},$$
where we write
$$(i',j')\leq (i,j)$$
if $i'\leq i$ and $j'\leq j$.
In the
present paper, we will consider $\Z\oplus\Z$-filtered
$\Z[U]$-complexes.  These come equipped with a chain map isomorphism
$U\colon C\ \longrightarrow C$
which carries $C\{(i,j)\}$ to $C\{(i-1,j-1)\}$.

Consider a subset 
$X\subset \Z\oplus\Z$ with the property that if $(i,j)\in X$, then
for any $(i',j')\leq (i,j)$, we also have that $(i',j')\in X$.If $C$
is any $\Z\oplus \Z$-filtered chain complex, we can form a subcomplex
$C\{X\}\subset C$ generated by $\bigoplus_{(i,j)\in X} C\{(i,j)\}$.

If $Y\subset \Z\oplus \Z$ a set with the property that for any
$(i,j)\in Y$, if $(i',j')\geq (i,j)$ we have that $(i',j')\in Y$. In
this case, we can endow $\bigoplus_{(i,j)\in Y} C\{(i,j)\}$ with the
structure of a quotient complex, which we will also denote by $C\{Y\}$.

If $K\subset S^3$ is a knot, we obtain an associated $\Z\oplus
\Z$-filtered chain complex $C$ with total homology isomorphic to
$\Z[U,U^{-1}]$.  The filtered chain homotopy type of this complex
$C=\CFKinf(S^3,K)$ is a knot invariant, \cite{Knots}, \cite{RasmussenThesis}.

The differential on $C$ induces also differential on each summand $C\{(i,j)\}$.
The homology group $H_*(C\{(0,s)\})$ is called the {\em knot Floer
homology group in filtration level $s$}, and it is denoted
$\HFKa(K,s)$.

The filtered chain homotopy type $\CFKinf(S^3,K)$ gives rise to some 
further algebraic structure.

Let $\Bp=C\{i\geq 0\}$. This is a model for $\CFp(S^3)$, and indeed,
so is $C\{j\geq 0\}$. There is a distinguished chain homotopy equivalence between
these two chain complexes.

We have also chain complexes 
$\Ap_s(K)=C\{\max(i,j-s)\geq 0\}$,
equipped with a pair of maps
\begin{eqnarray*}
\vertp_s\colon \Ap_s(K) \longrightarrow \Bp
&{\text{and}}&
\horp_s\colon \Ap_s(K) \longrightarrow \Bp,
\end{eqnarray*}
where the first is simply projection map (from $C\{\max(i,j-s)\geq 0\}$ to $C\{i\geq 0\}$), 
while the second is a
composite of the projection map $C\{\max(i,j-s)\geq 0\}$
to $C\{j\geq s\}$, followed by the identification with $C\{j\geq 0\}$
(induced by multiplication by $U^s$), followed by the chain homotopy
equivalence of this with $\Bp$. These maps are the data necessary
to calculate the Heegaard Floer homology of arbitrary Dehn fillings of $S^3$
along $K$.

Recall that the Heegaard Floer homology of $Y$
admits a direct sum splitting indexed by the set of $\SpinC$
structures over $Y$, which in turn is an affine space for $H^2(Y;\Z)$.
In particular, if $K\subset S^3$, 
then there is a splitting
$$\HFp(S^3_{p/q}(K))\cong \bigoplus_{i\in\Zmod{p}} \HFp(S^3_{p/q}(K),i).$$
Fix an integer $i$, and  consider the chain complexes
\begin{eqnarray*}
\BigAp_i=\bigoplus_{s\in\Z}(s,\Ap_{\lfloor \frac{i+ps}{q}\rfloor}(K))
&{\text{and}}&
\BigBp_i=\bigoplus_{s\in\Z}(s,\Bp),
\end{eqnarray*}
where here $\lfloor x\rfloor$ denotes the greatest integer smaller
than or equal to $x$.
We view the above chain homomorphisms $\vertp$ and $\horp$ as maps
\begin{eqnarray*}
\vertp\colon (s,\Ap_{\lfloor \frac{i+ps}{q}\rfloor}(K))\longrightarrow
(s,\Bp)
&{\text{and}}&
\horp\colon 
(s,\Ap_{\lfloor \frac{i+ps}{q}\rfloor}(K))
\longrightarrow
(s+1,\Bp).
\end{eqnarray*}
Adding these up, we obtain a chain map
$$\Dp_{i,p/q}\colon \BigAp_i \longrightarrow \BigBp_i;$$
i.e.
$$\Dp_{i,p/q} \{(s,a_s)\}_{s\in\Z}
=\{(s,b_s)\}_{s\in\Z},$$
where here
$$b_{s}
=\vertp_{\lfloor \frac{i+ps}{q}\rfloor}(a_s)
+\horp_{\lfloor \frac{i+p(s-1)}{q}\rfloor}(a_{s-1}).$$

Let $\Xp_{i,p/q}$ denote the mapping cone of $\Dp_{i,p/q}$.  Note that
$\Xp_{i,p/q}$ depends on $i$ only through its congruence class 
modulo $p$.  Note also that $\Ap_s$ and $\Bp_s$ are relatively
$\Z$-graded, and the homomorphisms $\vertp_s$ and $\horp_s$ respect
this relative grading. The mapping cone $\Xp_i$ can be endowed with a
relative grading, with the convention that $\Dp_{i,p/q}$ drops the
grading by one.

The following is proved (in somewhat more generality) 
in Theorem~\ref{RatSurg:thm:RationalSurgeries}
of~\cite{RatSurg}:

\begin{theorem}
\label{thm:SurgeryFormula}
  Let $K\subset S^3$ be a knot, and let $p, q$ be a pair of relatively
  prime integers. Then, there is an identification $\sigma\colon
  \Zmod{p}\longrightarrow \SpinC(S^3_{p/q}(K))$ such that for each
  $i\in \Zmod{p}$, there is a relatively graded isomorphism of groups
  $$\Phi_{K,i}\colon
  H_*(\Xp_{i,p/q}(K))\stackrel{\cong}\longrightarrow
  \HFp(S^3_{p/q}(K),\sigma(i)).$$ Indeed, there is a uniquely
  specified absolute grading on the subcomplex
  $\BigBp_i\subset\Xp_{i,p/q}(K)$ (which is independent of $K$) for
  which the map $\Phi_{O,i}$ is an isomorphism (where here $O$ is the
  unknot). With the corresponding induced grading on $\Xp_{i,p/q}(K)$,
  $\Phi_K$ becomes an absolutely graded isomorphism.
\end{theorem}

It is useful to note that there is a conjugation on knot Floer homology
related to the conjugation invariance on closed manifolds, cf. Equation~\eqref{eq:Conjugation}.
In the form which we need it, this is an isomorphism
$$\Psi\colon H_*(\Ap_s)\stackrel{\cong}{\longrightarrow} H_*(\Ap_{-s})$$
for all integers $s$. Indeed, under this isomorphism, we have a homotopy-commutative diagram commutative diagram
$$\begin{CD}
        \Ap_s @>{\Psi}>> \Ap_{-s}\\
        @V{\vertp_s}VV @VV{\horp_{-s}}V \\
        \Bp @>{=}>> \Bp
\end{CD}
$$

It is also useful to note that, although $\Xp_{i,p/q}(K)$ is a very
large chain complex, if we are interested in the homology in degrees
less than or equal to some fixed quantity $d$, then this is contained
in much smaller chain complex. More precisely, suppose that $p,q>0$.
Then, since $\vertp_s$ is an isomorphism for all sufficiently large
$s$ and $\horp_s$ is an isomorphism for all sufficiently small $s$, we
can consider the subsets
$$\BigAp_{i,\leq \sigma} = \bigoplus_{\{s\in\Z||s|\leq \sigma\}} (s,\Ap_{\lfloor \frac{i+ps}{q}\rfloor}(K)) \subset \BigAp_{i}$$
$$\BigBp_{i\leq \sigma} = \bigoplus_{\{s\in\Z|-\sigma<s\leq \sigma\}}
(s,\Bp)\subset \BigBp_{i}.$$
The map $\Dp_{i,p/q}$ induces a map from
$\BigAp_{i,\leq \sigma}$ to $\BigBp_{i,\leq \sigma}$, whose mapping
cone, denoted $\Xp_{i,p/q,\leq \sigma}$, is a quotient complex of
$\Xp_{i,p/q,\leq \sigma}$. Now, the homology of the chain complex
$\Xp_{i,p/q}(K)$ in degrees less than or equal to $d$ agrees with the
homology of its quotient complex $\Xp_{i,p/q,\leq \sigma}(K)$ for some
$\sigma$ depending on $d$.

\subsection{Examples}
\label{subsec:Examples}

For $K$ the unknot, $C$ has a single generator $a$ as a
$Z[U,U^{-1}]$-module, which is supported in filtration level $(0,0)$
and grading zero.  The differentials are trivial.

We let $T_\ell$ denote the left-handed trefoil,
$T_r$ denote the right-handed trefoil,
and $S$ denote the figure eight knot.

For $K=T_r$, we have that $C$ has three generators as a 
$\Z[U,U^{-1}]$-module,  $a$, $b$, $c$,  in filtration levels $(-1,0)$
$(0,0)$, and $(0,-1)$ respectively, with the differential
$D b = a+c, Da=Dc=0$. 

For $K=T_\ell$,  $C$ has three generators as a 
$\Z[U,U^{-1}]$-module,  $a$, $b$, $c$,  in filtration levels $(0,1)$
$(0,0)$, and $(1,0)$ respectively, with the differential
$D a = D c = b$ and $D b =0$.

Finally, for $K=S$, $C$ has five generators as a $\Z[U,U^{-1}]$ module,
$a,b,c,d,e$. Here $a$ is supported in filtration level $(1,1)$,
$b$ in $(0,1)$, $c$ in $(1,0)$, and $d$ and $e$ in $(0,0)$.
Differentials are given by 
$D a= b + c$, $Db = -Dc = d$, $Dd=De=0$.

\begin{figure}
\mbox{\vbox{\epsfbox{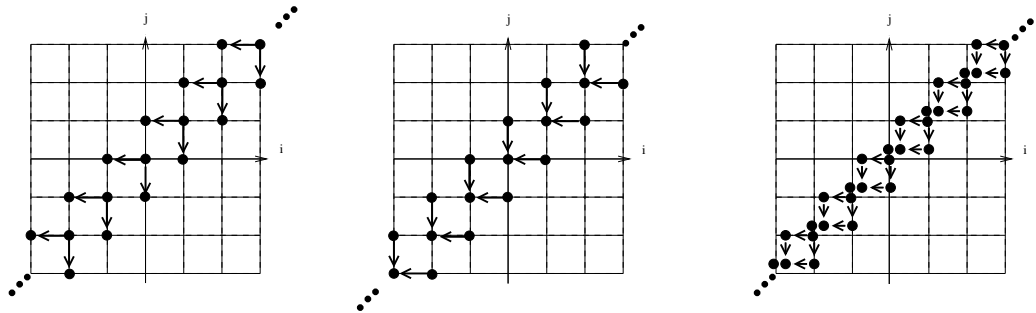}}}
\caption{\label{fig:MasterComplex}
{\bf{Filtered complexes for $T_r$, $T_\ell$, and $S$.}}
We have illustrated the $\Z\oplus\Z$-filtered complexes
associated to the three knots listed above. Dots represent generators,
and arrows represent differentials.}
\end{figure}

These answers are illustrated in Figure~\ref{fig:MasterComplex}.

These examples are all fairly simple to calculate; the calculations
can be made by finding appropriate Heegaard diagrams, or referring
to more general results from~\cite{RasmussenTwoBridge}, \cite{AltKnots}.

They have the following consequences for $\Ap_s$.

\begin{prop}
The groups $H_*(\Ap_s(K_0))$ and the homomorphisms $\vertp_s$ and $\horp_s$
for $K=T_r$, $T_\ell$, and $S$ are
determined  as follows.
For all $s\in \Z$, $H_*(\Ap_s(T_r))\cong \InjMod$; indeed,
for all $s>0$, the map induced by $v_s$ is an isomorphism, while for $s=0$, the map
$$v_0\colon H_*(\Ap_0(T_r))\cong\InjMod \longrightarrow \HFp(S^3)\cong\InjMod$$
is modeled on multiplication by $U$.

For all $s>0$, $H_*(\Ap_s(T_\ell))\cong \InjMod$, while $H_*(\Ap_0(T_\ell))\cong\InjMod_{(0)} \oplus \Z_{(0)}$, 
where the extra $\Z$ has grading zero. Moreover, the kernel of $v_0$ is one-dimensional.

For all $s>0$, $H_*(\Ap_s(T_\ell))\cong \InjMod$, while $H_*(\Ap_0(T_\ell))\cong\InjMod_{(0)} \oplus \Z_{(-1)}$, 
and the kernel of $v_0$ is one-dimensional.
\end{prop}

\begin{proof}
        These are straightforward consequences of the chain complexes described above.
\end{proof}

\subsection{Genus bounds and Ghiggini's theorem}

In~\cite{GenusBounds}, it is shown that if $K\subset S^3$ is a knot
with genus $g$, then 
$$\max\{s\big| H_*(C\{(0,s)\})=\HFKa(K,s)\neq 0\}=g.$$

It is conjectured that if $\HFKa(K,s)$ has rank one, then $K$ is a
fibered knot. In a recent paper, Ghiggini verified this conjecture 
for knots with genus one. Taken together, these results give the following:

\begin{theorem}
  Let $K$ be a knot with $\HFKa(K,s)=0$ for all $s>1$ and $\HFKa(K,1)$
  having rank one. Then, $K$ is a trefoil or the figure eight knot.
\end{theorem}

\section{Proof of Theorem~\ref{thm:Trefoil} and~\ref{thm:FigEight}}
\label{sec:Proof}

Let $K_0$ be a trefoil or the figure eight knot.  From a
graded isomorphism 
$$\HFp(S^3_{p/q}(K))\cong
\HFp(S^3_{p/q}(K_0)),$$ we would like to use the surgery formula to 
conclude an isomorphism between
the knot Floer homologies of $K$ and $K_0$. 
To this end, we find it useful to identify $\HFp(S^3_{p/q}(K_0))$.

\begin{prop}
        \label{prop:Calculations}
        Let $p$ and $q$ be relatively prime, positive integers.
        We have that 
        $$\rk \HFpRed(S^3_{p/q}(T_r))<q,$$ while
        $$\rk \HFpRed(S^3_{p/q}(T_\ell))=\rk \HFpRed(S^3_{p/q}(S))=q.$$
        In the case where $K_0=T_\ell$, $\HFpRed(S^3_{p/q}(T_\ell))$ is supported
        in even degree, while in $\HFp(S^3_{p/q}(S))$, it is supported in odd degree.
        Moreover, for $K_0=T_r$, $T_\ell$, and $S$, we have that
        $$\sum_{i}\rk \HFpRed(S^3_{p/q}(K_0),i) - \left(\frac{d(S^3_{p/q}(K_0),i)-d(S^3_{p/q}(O),i)}{2}\right)=q.$$
\end{prop}

\begin{proof}
  This is a straightforward application of the surgery formula
  (Theorem~\ref{thm:SurgeryFormula}) and the calculations from
  Subsection~\ref{subsec:Examples}.
\end{proof}

\begin{prop}
        \label{prop:NoBigGenus}
        If $K\subset S^3$ is a knot with the property that for some $s>0$, the map $H_*(\Ap_s(K)) \longrightarrow
        \HFp(S^3)$ induced by $\vertp$ is not an isomorphism, then 
        for any $p/q>0$, we have that
        $$\sum_{i}\rk
        \HFpRed(S^3_{p/q}(K),i) -
        \left(\frac{d(S^3_{p/q}(K),i)-d(S^3_{p/q}(O),i)}{2}\right)\geq
        2q.$$
\end{prop}

\begin{proof}
        Suppose that $\Hred(\Ap_s)\neq 0$. Then, we can find
        an element of $H_*(\Ap_s)$ with non-trivial image in
        $\Hred(\Ap_s)$ which is also contained in the kernel of the map
        $$\vertp_*\colon H_*(\Ap_s)\longrightarrow H_*(\Bp).$$
        From symmetry, 
        the same is true for 
        $$\horp\colon H_*(\Ap_{-s})\longrightarrow H_*(\Bp).$$
        
        Since for all $t$, the maps $\vertp_t$ and $\horp_t$ induce
        surjections on homology, it follows now from the surgery
        formula that the rank of $\HFpRed(S^3_{p/q}(K))$ is at least
        $2q$. 
        For example, suppose that is a cycle in $\Ap_s$
        representing a kernel of $\vertp_*([\xi_s])$, we claim that
        $(s,\xi_s)$ can be completed to a cycle class in $\Xp_{i,p/q}$
        by adding terms of the form $(t,\xi_t)\in \BigAp_{i}$ and
        $(t,\eta_t)\in\BigBp_{i}$ with $t>s$. For example, since
        $\vertp_{s+1}$ induces a surjection on homology, we can find a
        cycle $\xi_{s+1}\in\Ap_{\lfloor \frac{i+p(s+1)}{q}\rfloor}(K)$
        and a chain $\eta_{s+1}\in\Bp$ with
        $$\horp_{s}(\xi_s)=\vertp_{s+1}(\xi_{s+1})+\partial
        \eta_{s+1}.$$
        Proceeding inductively, we end up completing the
        initial cycle $\xi_s$ with a desired sequence of elements
        $\xi_t\in (t,\Ap_{\lfloor \frac{i+p(s+1)}{q}\rfloor}(K))$ and
        $\eta_t\in (t,\Bp)$ with $t\geq s$, so that
        $$\horp_{t}(\xi_t)=\vertp_{t+1}(\xi_{t+1})+\partial \eta_{t+1}.$$
        Moreover, for
        degree reasons, we can assume that the elements $\xi_t$ and
        $\eta_t$ vanish for sufficiently large $t$. Thus, the sum of
        these elements can be viewed as a homology class in
        $H_*(\Xp_{i,p/q})$ whose projection to
        $H_*(\Ap_{\lfloor\frac{i+ps}{q}\rfloor})$ is the given kernel
        element $[\xi_s]$ initial kernel element. This proves that if
        for $s>0$, $\Hred(\Ap_s)\neq 0$, we construct elements in
        $\Xp_{i,p/q}$, one for  for each time the chain complex
        $\Ap_s$ with appears ind $\Xp_{i,p/q}$, representing a homology class
        in $H_*(\Xp_{i,p/q})$ whose projection to
        $H_*(\Ap_s)$ has non-trivial image in $\Hred(\Ap_s)$. An
        analogous argument applies for the case where $s<0$.
        In view of this argument, the surgery formula now 
        guarantees that $\HFpRed(S^3_{p/q}(K))$ has rank at least $2q$.

        Suppose now that $\Hred(\Ap_s)=0$, but still the map induced
        by $\vertp_s$ is not an isomorphism.  It follows at once that
        the rank of its kernel $n$ is positive. Indeed, for all $0\leq
        t \leq s$, we see that $\vertp_t\Ap_t \longrightarrow \Bp$
        factors through the natural projection from $\Ap_t$ to
        $\Ap_s$.  It follows easily that the kernel of $\vertp_t
        \colon H_*(\Ap_t)\longrightarrow H_*(\Bp)$ for each such $t$
        has dimension at least $m$.

        Consider a $\SpinC$ structure $i\in \Zmod{p}$ modeled on  $\Xp_{p/q,i}(K)$.
        Let $k$ denote the number of copies of $\Ap_t$ with $|t|\leq s$ which appear in this model,
        and suppose that $k>0$. Then, it follows readily that
        $-(d(S^3_{p/q}(K),i)-d(S^3_{p/q}(O),i))\geq 2n$; and also that 
        $\rk \HFred(S^3_{p/q}(K)) \geq n(k-1)$.

        Thus, we see that
        $$\rk \HFpRed(S^3_{p/q}(K))-\left(\frac{d(S^3_{p/q}(K),i)-d(S^3_{p/q}(O),i)}{2}\right)\geq
        q(2|s|+1)n\geq 2q$$
\end{proof}

\begin{lemma}
        \label{lemma:KnotHomology}
        Let $K\subset S^3$ be a knot with genus $g$. Then, 
        there is a short exact sequence
        $$\begin{CD}
        0@>>> \HFKa(K,g)@>>>H_*(\Ap_{g-1})@>{\vertp_{g-1}}>> H_*(\Bp)@>>> 0
        \end{CD}$$
\end{lemma}

\begin{proof}
        There is an obvious short exact sequence
$$
        \begin{CD}
        0@>>>C\{(-1,g-1)\}@>>> \Ap_{g-1}@>{\vertp_{g-1}}>> \Bp @>>>0
        \end{CD}
$$
        inducing a long exact sequence in homology. On the other hand, the map
        on homology $\vertp_{g-1}$ is surjective for simple algebraic reasons.
        (as it is an isomorphism in all sufficiently large degrees, it is $U$-equivariant,
        and the automorphism of $H_*(\Bp)$ induced by $U$ is surjective.)
	Finally, note that $H_*(C\{(-1,g-1)\})\cong \HFKa(K,g)$.
\end{proof}

\noindent{\bf{Proof of Theorems~\ref{thm:Trefoil} and \ref{thm:FigEight}.}}
By reflecting the knot if necessary, we can assume that $p/q>0$. 
Assume that $K$ is a knot with
$S^3_{p/q}(K)\cong S^3_{p/q}(K_0)$, with $K_0\in\{T_r,T_\ell,S\}$. 
Combining Propositions~\ref{prop:Calculations} and \ref{prop:NoBigGenus},
we conclude that for all $s>0$,  $\vertp_s\colon \Ap\longrightarrow \Bp$ induces
an isomorphism on homology. From Lemma~\ref{lemma:KnotHomology}, we conclude that
$\HFKa(K,s)=0$ for all $s>1$. Note that this already proves that the genus of $K$ is one, 
cf.~\cite{GenusBounds}. 

Now, we claim that 
$$q\cm \rk \Ker\left( \vertp\colon H_*(\Ap_0(K))\longrightarrow
  H_*(\Bp)\right) = \rk \HFpRed(K) -\sum_i
  \left(\frac{d(S^3_{p/q}(K),i)-d(S^3_{p/q}(O),i)}{2}\right).$$

In view of Proposition~\ref{prop:Calculations}, we conclude that 
$\Ker\left( \vertp\colon H_*(\Ap_0(K))\longrightarrow H_*(\Bp)\right)$ has rank one. Thus,
by Lemma~\ref{lemma:KnotHomology}, we see that $\HFKa(K,1)$ has rank one. By Ghiggini's
theorem, it follows that $K$ is either the figure eight knot or the trefoil. 

Another look at the Floer homology groups $S^3_{p/q}(K_0)$ as stated in Proposition~\ref{prop:Calculations}
then allows one to conclude that $K=K_0$.
\qed

\vskip.2cm
\noindent{\bf{Proof of Corollary~\ref{cor:S237}.}}  Note that
$\Sigma(2,3,7)$ cannot be realized as $1/n$ surgery on any knot in
$S^3$. This follows from the surgery formula for Casson's invariant $\lambda$,
together with the fact that $|\lambda(\Sigma(2,3,7))|=1$.
The corollary is now a direct application of Theorems~\ref{thm:Trefoil} 
and \ref{thm:FigEight}.
\qed

\bibliographystyle{plain}
\bibliography{biblio}

\end{document}